\documentclass[oneside,10pt]{article}          

\usepackage[letterpaper]{geometry}	    
\usepackage{amsfonts,amsmath,latexsym,amssymb} 
\usepackage{theorem}                
\usepackage{mathrsfs,upref}         
\usepackage{mathptmx}		    

\newtheorem{thm}{Theorem}[section] 
\theoremstyle{definition} 
\newtheorem{defn}[thm]{Definition} 
\newtheorem{rem}[thm]{Remark} 
\newtheorem{exa}[thm]{Example} 
\numberwithin{equation}{section}

\title{Notes on Ordinary Conformal Differential Equations of order $\alpha$ ($0<\alpha\leq 1$)}
\author{Carlos E. Cadenas R.\\
	Dpto. de Matem\'aticas, FaCyT, Universidad de Carabobo.
	Valencia, Venezuela.\\
	Centro Multidisciplinario de Visualizaci\'on y C\'omputo Cient\'ifico,
	UC. Venezuela.}

\begin{document}

\maketitle

\thanks{$*$ Corresponding author. Email: ccadenas@uc.edu.ve}

\begin{abstract}
		These notes aim to provide a classical approach to solving some conformable differential equations based on prior knowledge of how to solve ordinary differential equations. That is, using the methods of separation of variables, homogeneous equations, linear, Bernoulli and exact. Representative examples are presented in all cases. Emphasis is placed on the new definitions and notations.
\end{abstract}
\maketitle

\begin{section}{Introduction}
	
	Since the publication of the article by Khalil et al. \cite{Khalil2014} there have been many publications related to the conformal derivative. Some of them are related to its development in the area of calculus (e.g. \cite{Abdeljawad2015,Atangana2015}). Many others are related to various applications \cite{ABas2019, Godinho2025, Khalil2014b, Nikolova2023, SALIM2026}. There are also detractors, who do not consider it to belong to the set of fractional derivatives. Without entering into that discussion, these notes aim to offer a summary of the methods for solving conformable ordinary differential equations using an approach similar to that classically used to solve ordinary differential equations, particularly the one followed in basic science and engineering courses. To this end, it has been necessary to provide the author's own definitions.	
	
	With this objective in mind, the present document is developed in the following sections. Section 2 presents some definitions and properties related to the conformal derivative. Section 3 gives the notation and some basic definitions necessary to perform the desired development. Section 4 introduces the method of separation of variables for  $\alpha-$separable equations. Section five defines when a differential equation is $(n,\alpha)$-homogeneous and how to solve it. Sections six and seven develop the  $(\alpha)-$linear and $(\alpha)-$Bernoulli differential equations, respectively. Section eight then develops  $(\alpha)-$exact differential equations. Representative examples are presented in all cases. Finally, a section providing conclusions and recommendations is presented.
	\begin{rem}	
		in the rest of these notes, as established in the title, $0<\alpha\leq 1$ will be considered.
	\end{rem}
	
\end{section}


\begin{section}{Conformable derivative and $\alpha-$integral}
	
	We will proceed to define the conformal derivative, which is due to Khalil et al., see \cite{Khalil2014}.
	\begin{defn}
		Given $f:\left[0,\infty \right) \rightarrow \mathbb{R}$. Then the ‘‘conformable fractional derivative’’ of $f$ of order $\alpha$ is defined by:
		\begin{equation}\label{conf}
			\frac{df(x)}{dx^{(\alpha)}}=\lim_{\epsilon \to 0}\frac{f(x+\epsilon x^{1-\alpha})-f(x)}{\epsilon}
		\end{equation}
		for all $x>0, \alpha\in(0,1)$.
	\end{defn}	
	In the original definition they use $T_\alpha$ as the operator which they call the fractional derivative of order $\alpha$, which is given by
	(\ref{conf}). In these notes we will use $\frac{df(x)}{dx^{(\alpha)}}, (f)^{(\alpha)}$ or $f^{(\alpha)}(x)$ indistinctly.
	
	Below we mention some of the properties of the conformal derivative that will be used as the document progresses. The following theorem can be found in \cite{Khalil2014}.
	\begin{thm}\cite{Khalil2014}
		Let $\alpha \in(0,1]$ and $f,g$ be $\alpha$-differentiable at a point $x>0$. Then
		\begin{enumerate}
			\item $(af+bg)^{(\alpha)}=a(f)^{(\alpha)}+b(g)^{(\alpha)}$, for all $a,b\in \mathbb{R}$.
			\item $(x^n)^{(\alpha)}=nx^{n-\alpha}$, for all $n\in \mathbb{R}$.
			\item $(C)^{(\alpha)}=0$, for all constant functions $f(x)=C$.
			\item $(fg)^{(\alpha)}=f(g)^{(\alpha)}+g(f)^{(\alpha)}$
			\item $\left(\frac{f}{g}\right)^{(\alpha)}=\frac{g(f)^{(\alpha)}-f(g)^{(\alpha)}}{g^2}$.
			\item If $f$ is differentiable, then $(f)^{(\alpha)}(x)=x^{1-\alpha}\frac{df}{dx}(x)$.	
		\end{enumerate}
	\end{thm}
	
	On the other hand, the chain rule for the conformable derivative is given by:
	
	\begin{thm}\cite{Atangana2015}
		Let $f$ and $g$ be two function differentiable, such that $g$ is differentiable at any $t$, and $f$ is differentiable at any $g(t)$, and then the conformable derivative obeys the Chain rule, meaning
	\end{thm}
	\begin{equation*}
		(f\circ g(x))^{(\alpha)}= x^{1-\alpha}g(x)^{1-\alpha}\frac{dg(x)}{dx} f^{(\alpha)}(t)\left|_{t=g(x)}\right. 
	\end{equation*}
	
	To determine whether a differential equation is exact, it is necessary to use Clairaut's theorem, also known as Schwarz's theorem or Young's theorem. Therefore, the definition of a conformable partial derivative and Clairaut's theorem for partial derivatives of fractional orders \cite{Atangana2015} are presented below.
	\begin{defn} \cite{Atangana2015}
		Let $f$ be a function with $m$ variables $x_1,\cdots, x_m$, and the conformable partial derivative of $f$ of order $0<\alpha\leq 1$ in $x_i$ is defined as follow
		\begin{equation*}
			\frac{\partial f(x_1,\cdots,x_i,\cdots, x_m)}{\partial x_i^{(\alpha)}}=\lim_{\epsilon \to 0}\frac{f(x_1,\cdots,x_i+\epsilon x_i^{1-\alpha},\cdots,x_m)-f(x_1,\cdots,x_i,\cdots, x_m)}{\epsilon}
		\end{equation*}
	\end{defn}
	
	\begin{thm}\cite{Atangana2015}
		Assume that $f(x,y)$ is function for which $\partial_x^{(\alpha)}\left[\partial_y^{(\beta)}f(x,y) \right]$ and $\partial_y^{(\beta)}\left[\partial_x^{(\alpha)}f(x,y) \right]$ exist and are continuous over the domain  $D\subset \mathbb{R}^2$ then
		
		\begin{equation} \label{Clairaut}
			\partial_x^{(\alpha)}\left[\partial_y^{(\beta)}f(x,y) \right]=\partial_y^{(\beta)}\left[\partial_x^{(\alpha)}f(x,y) \right]
		\end{equation}
	\end{thm}
	
	The $\alpha$-fractional integral of a function $f$ starting from $a\geq0$ is defined in \cite{Khalil2014} as follows:
	\begin{equation*}
		I_\alpha^a(f)(t)=\int_a^t \frac{f(x)}{x^{1-\alpha}}dx^{(\alpha)}
	\end{equation*}
	Therefore it is easy to obtain that $ \left( I_\alpha^a(f)(t)\right)^{(\alpha)} =f(t)$. This is presented in Theorem 3.1 of \cite{Khalil2014}.
	
	
	
\end{section}


\begin{section}{Preliminaries}
	
	From property 6 of theorem 1 it is clear that the relation between $dx$ and $dx^{(\alpha)}$ is given by
	
	\begin{equation}\label{dx}
		dx=x^{1-\alpha}dx^{(\alpha)}
	\end{equation}
	
	Let $f:\Omega \rightarrow \textbf{R}$ be differentiable in $(a,b)$, then we define the differential of $f$ as follows:
	\begin{equation}\label{difer}
		df=f^{(\alpha)}(x)dx^{(\alpha)}
	\end{equation}
	where $\alpha$ is the order of the fractional derivative with $0<\alpha\leq 1$. Note that if $\alpha=1$, the classical definition of the differential of a function holds.
	
	Now we will define the indefinite conformable integral.
	
	\begin{defn}
		Given a function $f$ an conformable fractional antiderivative of order $\alpha$ of $f$ is any function $F$ such that $F^{(\alpha)}(x)=f(x)$. 
		If $F$ is any conformable fractional antiderivative of order $\alpha$ then the conformable fractional antiderivative of order $\alpha$ de $f$ is called an  conformable fractional indefinite integral of order $\alpha$  and denoted,
		\begin{equation}\label{Int_frac}
			I^{(\alpha)}(f(x))=\int f(x)dx^{(\alpha)}=\int \frac{f(x)}{x^{1-\alpha}}dx=F(x)+C
		\end{equation}
	\end{defn}

\end{section}



Below are presented in Table 1 four conformable indefinite integrals that will be used in this document.

\begin{table}[htb]
	\begin{center}
		\begin{tabular}{| c | c | c |  c | c |}
			\hline
			$f(x)$ & $\int f(x)dx^{(\alpha)}$ & & $f(x)$ & $\int f(x)dx^{(\alpha)}$ \\ \hline		
			1 & $\frac{x^{\alpha}}{\alpha} + C$ & & $x^n, n\ne -\alpha$ &  $\frac{x^{n+\alpha}}{n+\alpha}+C$ \\ \hline
			$\frac{1}{x^\alpha}$ & $\ln\left| x\right| + C$ & &  $x^{1-\alpha}\frac{g'(x)}{g(x)}$ &  $\ln\left| g(x)\right| +C$ \\ \hline
		\end{tabular}
		\caption{
			Indefinite conformable integrals of classical functions.}
		\label{tab:Int_conf}
	\end{center}
\end{table}
Remark: It is clear that a more extensive table with conformable indefinite integrals could be presented. This was not done because its generation can be done with very basic knowledge of integral calculus.


In \cite{Abdeljawad2015} the author presents a theorem (Theorem 3.1) that establishes the integration by parts formula for the case of conformable definite integrals. A similar theorem is presented below for the case of conformable indefinite integrals.
\begin{thm}\label{Por parte} 
	Let $f$ and $g$ be two functions such that $fg$ is differentiable. Then
	\begin{equation}\label{Por partes}
		I^{(\alpha)}(f(x)g^{(\alpha)}(x))=f(x)g(x)-I^{(\alpha)}(g(x)f^{(\alpha)}(x))
	\end{equation}
\end{thm}

Remark: To proff this theorem, simply use the conformable derivative of $fg$ and then integrate it, just as it is customary to develop the classical indefinite integration by parts formula.

\textbf{Exercise:} using the formula (\ref{Por partes}) for indefinite conformable integration by parts, the reader is left to verify the following expression:
\begin{equation}\label{Ejercicio}
	\int x^me^{rx^{\alpha}/\alpha}dx^{(\alpha)}=\frac{x^m}{r}e^{rx^{\alpha}/\alpha}-\frac{m}{r}\int x^{m-\alpha}e^{rx^{\alpha}/\alpha}dx^{(\alpha)}
\end{equation}


\begin{section}{$(\alpha)$-separable differential equations}
	Let the fractional differential equation be
	\begin{equation}\label{EDF}
		M(x,y)dx^{(\alpha)}+N(x,y)dy=0
	\end{equation}
	It is said to be of $(\alpha)$-separable variables if by means of an integrating factor it can be transformed into the differential equation
	$$x^{1-\alpha}F(x)dx^{(\alpha)}+G(y)dy=0$$
	so that to solve it the differential equation must be integrated, obtaining:
	$$\int x^{1-\alpha}F(x)dx^{(\alpha)}+\int G(y)dy=C$$
	which due to the application of property (\ref{dx}) is equivalent to
	$$\int F(x)dx+\int G(y)dy=C$$
	
	Remark: To solve the conformable ordinary differential equation, it is not necessary to transform the conformable indefinite integral into a classical indefinite integral. However, both options have been given here.

	\begin{exa}  Solve the fractional differential equation:
		$$\frac{dy}{dx^{(\alpha)}}=\frac{x^{1-\alpha}e^x}{y\left( 1+e^x\right) }$$
		
		To solve this conformable differential equation, multiply by $ydx^{(\alpha)}$ and apply the definition (\ref{difer}), thus obtaining
		\begin{equation}\label{Ej1}
			\frac{x^{1-\alpha}e^x}{1+e^x}dx^{(\alpha)}-y dy=0
		\end{equation}
		Then, (\ref{dx}) is used and we obtain
		\begin{equation*}\label{Ej1sep}
			\frac{e^x}{1+e^x}dx-y dy=0
		\end{equation*}
		Which is solved in the book \cite{Kiseliov1968} and whose solution is $y^2-2\ln(1+e^x)=C$. Note that in equation (\ref{Ej1}) indefinite integration can also be applied directly and the same solution would be obtained by (\ref{Int_frac}). That is
		\begin{equation*}\label{Ej1a}
			\int \frac{x^{1-\alpha}e^x}{(1+e^x)}dx^{(\alpha)}-\int y dy=0
		\end{equation*}
		obtaining in the same way $y^2-2\ln(1+e^x)=C$ as a solution.
	\end{exa}

	\begin{exa} Solve the fractional differential equation: (example 4.2 of \cite{Khalil2014}): 
		$$\frac{dy}{dx^{(\alpha)}}+y=0$$
		By separating the variables we obtain
		$$\frac{dy}{y}+dx^{(\alpha)}=0$$
		Which when integrated indefinitely has the following solution:
		$$\ln |y|+\frac{x^\alpha}{\alpha}=C$$
		Or explicitly
		
		$$y=Ae^{-\frac{x^\alpha}{\alpha}}$$
	\end{exa}

	\begin{exa} Solve the conformable differential equation
		$$\frac{dy}{dx^{(\alpha)}}+x^\beta y=0,\quad \beta\in \mathbf{R}-\left\lbrace -\alpha\right\rbrace $$
		By separating the variables we obtain
		$$\frac{dy}{y}+x^\beta dx^{(\alpha)}=0$$
		That by integrating indefinitely and clearing the variable $y$ we have the solution:
		$$y=Ae^{-\frac{x^{\alpha+\beta}}{\alpha+\beta}}$$
		
		NOTE: If in this example $\beta=-\alpha$, the solution to the fractional differential equation is $xy=C$.
	\end{exa}

	\begin{subsection}{A particular case}
		Let the differential equation be
		\begin{equation}\label{CP1}
			\frac{dy}{dx^{(\alpha)}}=f(ax^{\alpha}+by+c)
		\end{equation}
		By means of a change of variables $z=ax^{\alpha}+by+c$ it can be reduced to a differential equation in separate variables. Then the derivative of $z$ of order $\alpha$ is given by
		\begin{equation}\label{DerCV1}
			\frac{dz}{dx^{(\alpha)}}=a{\alpha}+b\frac{dy}{dx^{(\alpha)}}
		\end{equation}
		Substituting (\ref{CP1}) into equation (\ref{DerCV1}) gives
		\begin{equation*}
			\frac{dz}{dx^{(\alpha)}}=a{\alpha}+bf(z)
		\end{equation*}
		separating variable
		
		\begin{equation*}
			\frac{dz}{a{\alpha}+bf(z)}=dx^{(\alpha)}
		\end{equation*}
		that by integrating indefinitely
		\begin{equation*}
			\int\frac{dz}{a{\alpha}+bf(z)}=\frac{x^\alpha}{\alpha}+C
		\end{equation*}
		After solving the integral of the left-hand side, we proceed to return the change and thus obtain the solution to the differential equation (\ref{CP1}).
		
	\end{subsection}
	
	\begin{exa} Solve the equation
		\begin{equation}\label{CP2}
			\frac{dy}{dx^{(\alpha)}}=\left( x^{\alpha}+y\right) ^2
		\end{equation}
		By using the change of variables $z=x^{\alpha}+y$ and differentiating $\frac{dz}{dx^{(\alpha)}}=\alpha+\frac{dy}{dx^{(\alpha)}}$. By using the equation (\ref{CP2})
		\begin{equation*}
			\frac{dz}{dx^{(\alpha)}}=\alpha+z^2
		\end{equation*}
		Separating variables
		
		\begin{equation}\label{Ecua}
			\frac{dz}{\alpha+z^2}=dx^{(\alpha)}
		\end{equation}
		By applying the indefinite integral to both sides of the differential equation (\ref{Ecua}) and solving said integrals, we obtain
		
		\begin{equation*}
			\sqrt{\alpha}\arctan\left( \frac{z}{\sqrt{\alpha}}\right) =\frac{x^{\alpha}}{\alpha}+C_1
		\end{equation*}
		After returning the change and solving for $y$ we have:
		
		\begin{equation*}
			y=\sqrt{\alpha}\tan\left( \frac{x^{\alpha}}{\alpha^{3/2}}+C\right)-x^{\alpha}
		\end{equation*}
		This example is an adaptation of a similar one presented in \cite{Krasnov1990}
		
	\end{exa}
	
\end{section}

\begin{section}{Differential equations $(n,\alpha)$-homogeneous}
	
	The function $M$ is defined to be homogeneous of order $n+1-\alpha$ if $M(\lambda x,\lambda y)=\lambda^{n+1-\alpha}M(x,y)$. Similarly, $N$ is homogeneous of order $n$ if $N(\lambda x,\lambda y)=\lambda^{n}N(x,y)$. Now, the fractional differential equation given in (\ref{EDF}) is said to be homogeneous of order $(n,\alpha)$ if $M$ and $N$ are homogeneous of order $n+1-\alpha$ and $n$, respectively. In this case, we can proceed using the change of variable $y=xu$, so using the derivative rule of a product, we have	
	
	$$\frac{dy}{dx^{(\alpha)}}=\frac{dx}{dx^{(\alpha)}}u+x\frac{du}{dx^{(\alpha)}}$$
	or what is the same
	$$\frac{dy}{dx^{(\alpha)}}=x^{1-\alpha}u+x\frac{du}{dx^{(\alpha)}}$$
	or in the differential form
	
	$$dy=x^{1-\alpha}u dx^{(\alpha)}+x du$$	
	Replacing the variable change mentioned in (\ref{EDF}):
	$$M(x,xu)dx^{(\alpha)}+N(x,xu)(x^{1-\alpha}u dx^{(\alpha)}+x du)=0$$
	Assuming that the fractional differential equation given in (\ref{EDF}) is homogeneous of order $(n,\alpha)$ then the above equation becomes:
	$$x^{n+1-\alpha}M(1,u)dx^{(\alpha)}+x^{n}N(1,u)(x^{1-\alpha}u dx^{(\alpha)}+x du)=0$$
	After separating variables, considering $F(u)=M(1,u)$ and $G(u)=N(1,u)$ we have:
	$$\frac{dx^{(\alpha)}}{x}+\frac{G(u)du}{F(u)+uG(u)}=0$$
	By integrating this differential equation
	
	$$\frac{x^{\alpha-1}}{\alpha-1}+\int\frac{G(u)du}{F(u)+uG(u)}=C, \quad (\alpha\neq 1)$$
	
	If $\alpha= 1$, then
	$$\ln|x|+\int\frac{G(u)du}{F(u)+uG(u)}=C, \quad (\alpha\neq 1)$$
	After solving the indefinite integral, the change must be returned.

	\begin{exa}
		
		The differential equation $(2x+3y)dy=x^{1-\alpha}(x+y)dx^{(\alpha)}$ is $(1,\alpha)$-homogeneous which is equivalent to the ordinary differential equation
		$$\frac{dy}{dx}=\frac{x+y}{2x+3y}$$
		which appears in example 4.4 of \cite{Khalil2014}.
		
	\end{exa}

	\begin{subsection}{A particular case}
		Consider the fractional differential equation $(0,n)$-homogeneous
		$$\frac{dy}{dx^{(\alpha)}}=x^{1-\alpha}\psi(y/x)$$
		When using the change of variables $y=ux$ we have that
		$$\frac{dy}{dx^{(\alpha)}}=x^{1-\alpha}u+x\frac{du}{dx^{(\alpha)}}=x^{1-\alpha}\psi(u)$$
		Separating variables
		$$\frac{dx^{(\alpha)}}{x^\alpha}+\frac{du}{u-\psi(u)}=0$$
		By integrating indefinitely we have
		$$\ln|Cx|+\int\frac{du}{u-\psi(u)}=0$$
		
		\begin{exa}
			Solve the differential equation
			$$\frac{dy}{dx^{(\alpha)}}=x^{1-\alpha}\left(\frac{y}{x}+\frac{x}{y} \right) $$
			In this case $\psi(u)=u+1/u$, so solving this fractional differential equation would give
			$$\ln|Cx|+\frac{u^2}{2}=0$$
			and when returning the change you have to $y^2+x^2\ln{|Cx^2|}=0$
			
		\end{exa}
		
	\end{subsection}
	
	\begin{subsection}{Convertible to $(1,\alpha)$-homogeneous}
		
		Let $h$ and $k$ be the point of intersection between the lines $a_1x+b_1y+c_1=0$ and $a_2x+b_2y+c_2=0$, the equation
		$$(x-h)^{1-\alpha}(a_1x+b_1y+c_1)dx^{(\alpha)}+(a_2x+b_2y+c_2)dy=0$$
		It is converted to $(1,\alpha)$-homogeneous using the changes of variables $x=u+h$ and $y=v+k$
		
		\begin{exa}
			Solve the differential equation
			$$\frac{dy}{dx^{(\alpha)}}=-(x-1)^{1-\alpha}\frac{2x-3y+4}{3x-2y+1}. $$
			Using the changes of variables $x=u+1$ and $y=v+2$ the following ODE is obtained
			$$u^{1-\alpha}(u+v)du^{(\alpha)}+(u-v)dy=0$$
			This ODE is $(1,\alpha)$-homogeneous, so by making the change of variable $v=zu$ and considering that $dv=u^{1-\alpha}zdu^{\alpha}+udz$ one has
			
			$$u^{1-\alpha}(2u-3zu)du^{(\alpha)}+(3u-2zu)(u^{1-\alpha}zdu^{\alpha}+udz)=0$$
			which when separating variables can be written as follows
			$$2\frac{du^{(\alpha)}}{u^{\alpha}}+\frac{3-2z}{1-z^2}dz=0$$
			whose solution is
			
			$$u^4(z+1)^5=C(z-1)$$
			returning the change in $z$ we have
			$$(v+u)^5=C(v-u)$$

			By returning the changes to the original variables, the solution is obtained
			$$(x+y-3)^5=C(y-x-1)$$
			This example is an adaptation of an example found in \cite{Kells1965}.
			
		\end{exa}
		
	\end{subsection}

\end{section}

\begin{section}{($\alpha$)-linear differential equations}
	Let the fractional differential equation ($\alpha$)-linear be:
	\begin{equation*}\label{EcLin}
		\frac{dy}{dx^{(\alpha)}}+P(x)y=Q(x)
	\end{equation*}	
	By multiplying it by the integrating factor $e^{I^{(\alpha)}(P(x))}$ is obtained
	\begin{equation*}\label{EcLin2}
		e^{I^{(\alpha)}(P(x))} dy+e^{I^{(\alpha)}(P(x))}P(x)y dx^{(\alpha)}=e^{I^{(\alpha)}(P(x))} Q(x) dx^{(\alpha)}
	\end{equation*}
	The left-hand side is the differential of a product, so
	$$d\left( e^{I^{(\alpha)}(P(x))} y\right) =e^{I^{(\alpha)}(P(x))} Q(x) dx^{(\alpha)}$$
	and by integrating indefinitely
	$$ e^{I^{(\alpha)}(P(x))} y =\int e^{I^{(\alpha)}(P(x))} Q(x) dx^{(\alpha)}+C$$
	when obtaining the variable $y$ we have
	$$ y =e^{-I^{(\alpha)}(P(x))}\left( I^{(\alpha)} \left( e^{I^{(\alpha)}(P(x))} Q(x)\right) +C\right) $$
	The next two examples are adaptations of examples 2, 3 and 4 of \cite{Kareem2017}.
	
	\begin{exa}
		
		Solve the differential equation
		\begin{equation*}\label{Ejx}
			\frac{dy}{dx^{(\alpha)}}+y=x^{m+1-\alpha}+(m+1-\alpha)x^{m+1-2\alpha}
		\end{equation*}
		Using the integrating factor $e^{\frac{x^\alpha}{\alpha}}$, you have to
		$$d\left( e^{\frac{x^\alpha}{\alpha}} y\right) =x^{1-\alpha}\left( x^{m}+(m+1-\alpha)x^{m-\alpha}\right)  e^{\frac{x^\alpha}{\alpha}}dx^{(\alpha)} $$
		Integrating and calculating the variable $y$
		
		$$y=x^{m+1-\alpha}+Ae^{-\frac{x^\alpha}{\alpha}}$$
		
	\end{exa}

	\begin{exa}
		\begin{equation*}\label{Ejy}
			\frac{dy}{dx^{(\alpha)}}+x^\beta y=x^me^{-\frac{x^{\alpha+\beta}}{\alpha+\beta}};\quad (\alpha\ne \beta)
		\end{equation*}
		Using the integrating factor  $e^{\frac{x^{\alpha+\beta}}{\alpha+\beta}}$,  you have to
		
		$$d\left( e^{\frac{x^{\alpha+\beta}}{\alpha+\beta}} y\right) =x^m e^{\frac{x^\alpha}{\alpha}}dx^{(\alpha)} $$
		Integrating and calculating the variable $y$
		
		$$y=\left( \frac{x^{m+\alpha}}{m+\alpha}+C\right)e^{-\frac{x^{\alpha+\beta}}{\alpha+\beta}} $$

		If $\alpha= \beta=\frac{1}{2}$ y $m=1$, example 4 is obtained from
		\cite{Kareem2017}.
		
	\end{exa}
	
\end{section}

\begin{section}{ ($\alpha$)-Bernoulli differential equations}
	
	In this section, using the same strategy as for classical differential equations, we will transform the ($\alpha$)-Bernoulli equation into a ($\alpha$)-linear equation. A representative example will then be presented.
	
	Let the ($\alpha$)-Bernoulli differential equation be
	
	\begin{equation*}\label{EcBern}
		\frac{dy}{dx^{(\alpha)}}+P(x)y=Q(x)y^n
	\end{equation*}	
	
	By multiplying the ODE by $(1-n)y^{-n}$ and making the change of variable $z=y^{1-n}$, as $\frac{dz}{dx^{(\alpha)}}=(1-n)y^{-n}\frac{dy}{dx^{(\alpha)}}$ is obtained
	\begin{equation*}\label{EcBernz}
		\frac{dz}{dx^{(\alpha)}}+(1-n)P(x)z=(1-n)Q(x)
	\end{equation*}	
	which is  ($\alpha$)-linear in $z$, and has as solution
	$$ z=y^{1-n} =e^{-I^{(\alpha)}((1-n)P(x))}\left( I^{(\alpha)} \left( (1-n)e^{I^{(\alpha)}((1-n)P(x))} Q(x)\right) +C\right) $$

	\begin{exa}
		\begin{equation*}\label{Eja}		
			\frac{dy}{dx^{(\alpha)}}+y=\left(x^r-\frac{r}{n-1}x^{r-\alpha}\right)y^n\quad (n\ne 1)
		\end{equation*}
		Multiplying the ODE by $(1-n)y^{-n}$ and making the change of variable $z=y^{1-n}$, we obtain
		\begin{equation*}\label{Ejaz}
			\frac{dz}{dx^{(\alpha)}}+(1-n)z=(1-n)\left(x^r-\frac{r}{n-1}x^{r-\alpha}\right)
		\end{equation*}
		Using the integrating factor $e^{(1-n)\frac{x^\alpha}{\alpha}}$ we obtain
		\begin{equation*}\label{Ejaz1}
			d\left(e^{(1-n)\frac{x^\alpha}{\alpha}} z\right)=(1-n)\left(x^r-\frac{r}{n-1}x^{r-\alpha}\right)e^{(1-n)\frac{x^\alpha}{\alpha}}dx(\alpha)
		\end{equation*}
		
		By integrating and obtaining $z$	
		\begin{equation}\label{Ejaz2}
			z=e^{(n-1)\frac{x^\alpha}{\alpha}}\left((1-n)\int\left(x^r-\frac{r}{n-1}x^{r-\alpha}\right)e^{(1-n)\frac{x^\alpha}{\alpha}}dx^{(\alpha)}+C\right)
		\end{equation}
		
		From exercise of integration by parts using (\ref{Ejercicio}) with $r=1-n$ we have that
		\begin{equation*}\label{EjPPaz2}
			\int\left(x^m-\frac{m}{n-1}x^{m-\alpha}\right)e^{(1-n)\frac{x^\alpha}{\alpha}}dx^{(\alpha)}=\frac{x^m}{1-n}e^{(1-n)\frac{x^\alpha}{\alpha}}+C
		\end{equation*}
		So by replacing $z$ in (\ref{Ejaz2}) we have
		\begin{equation*}\label{Ejaz2Sol}
			y^{1-n}=e^{(n-1)\frac{x^\alpha}{\alpha}}\left((1-n)\frac{x^r}{1-n}e^{(1-n)\frac{x^\alpha}{\alpha}}+C\right)
		\end{equation*}
		If $r=n=2$ and $\alpha=\frac{1}{2}$ we have example 3 of \cite{Ilie2017}.
		
	\end{exa}
	
\end{section}

\begin{section}{{($\alpha$)-exact} differential equations }
	
	If $z=f(x,y)=C$, the differential is given by
	\begin{equation}\label{EcDif2Var}
		dz^{(\alpha)}=\partial_x^{(\alpha)}f(x,y)dx^{(\alpha)}+\partial_yf(x,y)dy=0
	\end{equation}
	By Clairaut's theorem, for $\beta=1$ in equation (\ref{Clairaut}), if
	\begin{equation*}\label{EqClaEX}
		\partial_x^{(\alpha)}\left[\partial_yf(x,y) \right]=\partial_y\left[\partial_x^{(\alpha)}f(x,y) \right]
	\end{equation*}
	We say that the equation (\ref{EcDif2Var}) is ($\alpha$)-exact.
	
	Then, to solve the differential equation (\ref{EDF}) expressed in the form
	\begin{equation*}\label{EDFExacta}
		\partial_x^{(\alpha)}f(x,y) dx^{(\alpha)}+\partial_yf(x,y)dy=0
	\end{equation*}
	It must be verified that:
	\begin{equation*}\label{EcCondExacta}
		\partial_y M(x,y)=\partial_x^{(\alpha)} N(x,y)dy
	\end{equation*}
	which would be the condition for the differential equation (\ref{EDF}) to be ($\alpha$)-exact. Therefore, this differential equation is ($\alpha$)-exact.
	So, to solve the exact ($\alpha$) equation, proceed as follows:

	Since $M(x,y)=\partial_x^{(\alpha)}f(x,y)$ then
	\begin{equation*}\label{Ecfexactax}
		f(x,y)= \int M(x,y) dx^{(\alpha)}+g(y)
	\end{equation*}
	then, as it is satisfied that
	\begin{equation*}\label{EcNexactax}
		N(x,y)=\partial_y f(x,y)
	\end{equation*}
	so 
	\begin{equation*}\label{EcNexactax}
		N(x,y)=\partial_y \left( \int M(x,y) dx^{(\alpha)}+g(y) \right) 
	\end{equation*}
	
	Solving for $g(y)$ gives the solution to the ordinary differential equation of order $\alpha$: $f(x,y)+C$; that is
	\begin{equation*}\label{Ecfexactax}
		f(x,y)= \int M(x,y) dx^{(\alpha)}+\int\left( N(x,y)-\partial_y\int M(x.y) dx^{(\alpha)}\right) dy=C
	\end{equation*}
	
	\begin{exa}
		\begin{equation*}\label{Ejw}
			x^{1-\alpha}e^{-y}dx^{(\alpha)}-\left(2y+xe^{-y} \right) dy=0
		\end{equation*}

		Identifying $M$ and $N$ we have that: $M(x,y)=x^{1-\alpha}e^{-y}$ and $N(x,y)=-\left(2y+xe^{-y} \right)$, so their respective derivatives are
		$$\partial_y M(x,y)=\partial_x^{(\alpha)} N(x,y)=-x^{1-\alpha}e^{-y}$$
		
		So the differential equation in this example is ($\alpha$)-exact. Solving according to the above, we have:
		$$\partial_x^{(\alpha)} f(x,y)=x^{1-\alpha}e^{-y}=M(x,y)$$
		calculating $f$ we have:
		$$f(x,y)=\int{x^{1-\alpha}e^{-y}dx^{(\alpha)}}=xe^{-y}+g(y)$$
		
		On the other hand, as
		$$N(x,y)=-\left(2y+xe^{-y} \right)=\partial_y f(x,y)$$
		then
		$$-\left(2y+xe^{-y} \right)=-xe^{-y}+g'(y)$$
		then $g'(y)=-2y$, so $g(y)=-y^2+C_1$. Thus the solution to the differential equation is
		$$f(x,y)=xe^{-y}-y^2=C$$
		
	\end{exa}
	
\end{section}

\begin{section}{Final Comments}
	In these notes, a somewhat novel approach to solving conformable ordinary differential equations was presented. This approach was recently used in a small part of an elective course for a Bachelor's degree in Mathematics at the Faculty of Science and Technology of the University of Carabobo. The idea was to present all the basic topics related to the subject in a single document. To achieve this, it was necessary to establish some novel ideas, such as establishing the differential equation in the form $M(x,y)dx^{(\alpha)}+N(x,y)dy=0$ and thus presenting novel definitions such as ($\alpha$)-separable differential equations, (n,$\alpha$)-homogeneous, ($\alpha$)-linear, ($\alpha$)-Bernoulli, and ($\alpha$)-exact. 
\end{section}



\begin{thebibliography}{20}
	
	
	\bibitem{Abdeljawad2015}
	{\sc T. Abdeljawad},
	{\it On conformable fractional calculus},
	Journal of computational and Applied Mathematics. {\bf 279}, (2015), 57--66.
	
	\bibitem{Atangana2015}
	{\sc A. Atangana, D. Baleanu and A. Alsaedi},
	{\it New properties of conformable derivative},
	Open Mathematics. {\bf 13}, 1 (2015), 889--898. https://doi.org/10.1515/math-2015-0081        
	
	\bibitem{ABas2019}
	{\sc E. Bas, B. Acay and R. Ozarslan},
	{\it The price adjustment equation with different types of conformable derivatives in market equilibrium},
	AIMS Math. {\bf 4}, 3 (2019), 805--820.
	
	\bibitem{Godinho2025}
	{\sc C. F. de Lima Godinho, C. M. Porto, M. C. Rodriguez,  and I. V. Vancea},
	{\it  Coherent States of Conformable Quantum Oscillator},
	(2025).
	
	\bibitem{Ilie2017}
	{\sc M. Ilie, J. Biazar and Z. Ayati},
	{\it General solution of Bernoulli and Riccati fractional differential equations based on conformable fractional derivative},
	Int. J. Appl. Math. Res. {\bf 6}, 2 (2017), 49--51.
	
	
	\bibitem{Kareem2017}
	{\sc A. M. Kareem},
	{\it Conformable fractional derivatives and it is applications for solving fractional differential equations},
	J. Math. {\bf 13}, 1 (2017), 81--87.
	
	
	\bibitem{Kells1965}
	{\sc L. Kells},
	{\it Elementary differential equations},
	Mc Graw-Hill (1965).        
	
	\bibitem{Khalil2014}
	{\sc  R. Khalil, et al},
	{\it A new definition of fractional derivative},
	Journal of Computational and Applied Mathematics. {\bf 264}, (2014), 65--70.
	
	\bibitem{Khalil2014b}
	{\sc R. Khalil and M. Abu-Hammad},
	{\it Conformable fractional heat differential equation},
	International Journal of Pure and Applied Mathematics,   {\bf 94}, 2 (2014), 215-221.
	
	\bibitem{Kiseliov1968}
	{\sc A. Kiseliov,  et al},
	{\it Problemas de ecuaciones diferenciales ordinarias},
	Editorial Mir Moscow. (1968).
	
	\bibitem{Krasnov1990}
	{\sc M. Krasnov, et al},
	{\it Mathematical analysis for engineers},
	Vol. 2. Mir Publishers Moscow. (1990).
	
	\bibitem{Nikolova2023}
	{\sc E. V. Nikolova and M. Chilikova–Lubomirova},
	{\it  On the Solitary Wave Solutions of the Conformable Time-Fractional Wu–Zhang System via Simple Equations Method (SEsM)}
	In Annual Meeting of the Bulgarian Section of SIAM, Cham: Springer Nature Switzerland. (2023) 164-174.
	
	
	\bibitem{SALIM2026}
	{\sc A. Salim, S. Krim, J. E. Lazreg and M. Benchohra},
	{\it A study on some conformable fractional implicit hybrid differential equations with delay},
	Kragujevac Journal of Mathematics, {\bf 50}, (3) (2026) 439-455.

	
\end{thebibliography}
\end{document}